\documentclass{article}
\usepackage{amssymb}
\usepackage{graphicx}
\usepackage{amsmath}

\input{tcilatex}
\begin{document}

\begin{center}
\textbf{\ Nonlocal fractional differential equat\i ons and applications\ \ }

\textbf{Veli Shakhmurov }

\ \ \ \ \ \ \ \ Istanbul Okan University, Department of Mechanical
Engineering, Akfirat, Tuzla 34959 Istanbul, Turkey, E-mail:
veli.sahmurov@okan.edu.tr;

\ \ \ \ \ \ \ \ \ \ \ \ \ \ 
\end{center}

\begin{quote}
\ \ \ \ \ \ \ \ \ \ \ \ \ \ \ \ \ \ \ \ \ \ \ \ \ \ \ \ \ \ \ \ \ \ \ \ \ \
\ \ \ \ \ \ \ \ \ \ \ \ \ \ \ 
\end{quote}

\begin{center}
\textbf{ABSTRACT}
\end{center}

\begin{quote}
\ \ \ \ \ \ \ \ \ \ \ \ \ \ \ \ \ 
\end{quote}

\ \ Boundary value problems for\ nonlocal fractional elliptic equations with
parameter in Banach spaces are studied. Uniform $L_{p}$-separability
properties and sharp resolvent estimates are obtained for elliptic equations
in terms of fractional derivatives. Particularly, it is proven that the
fractional ellipitic operator generated by these equations is sectorial and
also is a generator of an analytic semigroup. Moreover, maximal regularity
properties of\ nonlocal fract\i onal abstract parabolic equation are
established. As an application, the nonlocal anisotropic fract\i onal\
differential equations and the system of nonlocal fract\i onal\ differential
equations are studied.\ \ \ \ \ 

\begin{center}
\bigskip\ \ \textbf{AMS: 47GXX, 35JXX, 47FXX, 47DXX, 43AXX}
\end{center}


\begin{center}
\textbf{1. Introduction, notations and background }
\end{center}

In the last years, the maximal regularity properties of boundary value
problems (BVPs) for abstract differential equations (ADEs) have found many
applications in PDE and pseudo DE with applications in physics (see $\left[ 
\text{1, 2}\right] $, $\left[ \text{5, 6}\right] $, $\left[ \text{8, 9}%
\right] $, $\left[ \text{15-19}\right] $,$\left[ \text{22,23}\right] $ and
the references therein). ADEs have found many applications in fract\i onal
differential equations (FDEs), pseudo-differential equations (PsDE) and
PDEs. FDEs were treated e.g. in $\left[ \text{4, 7, 10-12, 14, 21}\right] $.
The regularity properties of FDEs have been studied e.g. in $\left[ \text{%
11, 12, 20}\right] $. The exsistence and uniqueness of solution to
fractional ADEs were studied e.g. in $\left[ \text{3, 11, 12}\right] $.
Regularity properties of nonlocal ADEs were investigated e.g. in $\left[
15-17\right] $. The main objective of the present paper is to discuss the  $%
L_{p}\left( \mathbb{R};H\right) $-maximal regularity of the fractional\ ADE
with parameter%
\begin{equation}
a\ast D^{\gamma }u+A\ast u+\lambda u=f\left( x\right) ,\text{ }x\in \mathbb{%
R=}\left( -\infty ,\infty \right) ,  \tag{1.1}
\end{equation}%
where $a$ is complex valued functions, $\lambda $ is a complex parameter, $%
A=A\left( x\right) $ is a linear operator functions in a Hilbert space $H$,
and $D^{\gamma }$is Riemann-Liouville type fractional derivatives of order $%
\gamma \in \left( \left. 1,2\right] ,\right. $ i.e. 
\begin{equation}
D^{\gamma }u=\frac{1}{\Gamma \left( 2-\gamma \right) }\frac{d^{2}}{dx^{2}}%
\dint\limits_{0}^{x}\frac{f\left( y\right) dy}{\left( x-y\right) ^{\gamma -1}%
},  \tag{1.2}
\end{equation}%
here $\Gamma \left( \gamma \right) $ is Gamma function for $\gamma >0$ (see
e.g. $\left[ \text{7, 10}\right] $) and the convolutions $a\ast D^{\gamma }u$%
, $\ A\ast u$ are defined in the distribution sense (see e.g. $\left[ \text{%
1, Section 3}\right] $).

For\ $\alpha _{i}\in \left[ 0,\infty \right) $ and $\mathbf{\alpha }=\left(
\alpha _{1},\alpha _{2},...,\alpha _{n}\right) $. Here, $D^{\alpha
}=D_{1}^{\alpha _{1}}D_{2}^{\alpha _{2}},.,D_{n}^{\alpha _{n}}$. Let $E$ be
Banach space. Here, $L_{p}\left( \Omega ;E\right) $ denotes the space of
strongly measurable $E$-valued functions that are defined on the measurable
subset $\Omega \subset \mathbb{R}^{n}$ with the norm given by

\begin{equation*}
\left\Vert f\right\Vert _{L_{p}\left( \Omega ;E\right) }=\left(
\int\limits_{\Omega }\left\Vert f\left( x\right) \right\Vert
_{E}^{p}dx\right) ^{\frac{1}{p}},\text{ }1\leq p<\infty \ .
\end{equation*}

\bigskip Let $S\left( \mathbb{R}^{n};E\right) $ denote the $E$-valued
Schwartz class, i.e., the space of all $E$-valued rapidly decreasing smooth
functions on $\mathbb{R}^{n}$ equipped with its usual topology generated by
seminorms.

A function $\Psi \in C\left( \mathbb{R}^{n};B\left( E_{1},E_{2}\right)
\right) $ is called a Fourier multiplier from $L_{p}\left( \mathbb{R}%
^{n};E_{1}\right) $\ to $L_{p}\left( \mathbb{R}^{n};E_{2}\right) $ if the map%
\begin{equation*}
u\rightarrow \Lambda u=F^{-1}\Psi \left( \xi \right) Fu,\text{ }u\in S\left( 
\mathbb{R}^{n};E_{1}\right)
\end{equation*}%
is well defined and extends to a bounded linear operator 
\begin{equation*}
\Lambda :\ L_{p}\left( \mathbb{R}^{n};E_{1}\right) \rightarrow \ L_{p}\left( 
\mathbb{R}^{n};E_{2}\right) .
\end{equation*}

Let $h\in R,$ $m\in \mathbb{N}$ and $e_{k}$, $k=1,2,...,n$ be the standard
unit vectors of $\mathbb{R}^{n},$ 
\begin{equation*}
\Delta _{k}\left( h\right) f\left( x\right) =f\left( x+he_{k}\right)
-f\left( x\right) .
\end{equation*}

Let $A=A\left( x\right) ,$ $x\in \mathbb{R}^{n}$ be closed linear operator
in $E$ with domain $D\left( A\right) $ independent of $x.$ The Fourier
transformation of $A\left( x\right) $ is a linear operator with the domain $%
D\left( A\right) $ defined as 
\begin{equation*}
\hat{A}\left( \xi \right) u\left( \varphi \right) =A\left( x\right) u\left( 
\hat{\varphi}\right) \text{ for }u\in S^{\prime }\left( \mathbb{R}%
^{n};E\left( A\right) \right) ,\text{ }\varphi \in S\left( \mathbb{R}%
^{n}\right) .
\end{equation*}%
(For details see e.g $\left[ \text{1, Section 3}\right] $). $A\left(
x\right) $ is differentiable if there is the limit

\begin{equation*}
\left( \frac{\partial A}{\partial x_{k}}\right) u=\lim_{h\rightarrow 0}\frac{%
\Delta _{k}\left( h\right) A\left( x\right) u}{h},\text{ }k=1,2,...n,\text{ }%
u\in D\left( A\right)
\end{equation*}%
in the sense of $E-$norm.

We prove that problem $\left( 1.1\right) $ has a maximal regular unique
solution $u\in $\ $W_{p}^{\gamma }\left( \mathbb{R};H\left( A\right)
,H\right) $ for $f\in $ $L_{p}\left( \mathbb{R};H\right) $ and the following
uniform coercive estimate holds 
\begin{equation}
\dsum\limits_{0\leq s\leq \gamma }\left\vert \lambda \right\vert ^{1-\frac{s%
}{\gamma }}\left\Vert D^{s}u\right\Vert _{L_{p}\left( \mathbb{R};H\right)
}+\left\Vert A\ast u\right\Vert _{L_{p}\left( \mathbb{R};H\right) }\leq
C\left\Vert f\right\Vert _{L_{p}\left( \mathbb{R};H\right) }.  \tag{1.3}
\end{equation}%
\ The estimate $\left( 1.3\right) $ implies that the operator $O$ generated
by $\left( 1.1\right) $ has a bounded inverse from $L_{p}\left( \mathbb{R}%
;H\right) $ into the space $W_{p}^{\gamma }\left( \mathbb{R};H\left(
A\right) ,H\right) ,$ which will be defined subsequently. Particularly, from
the estimate $\left( 1.3\right) ,$ we obtain that the operator $O$ is
uniformly positive in $L_{p}\left( \mathbb{R};H\right) .$ By using this
property we prove the well posedness of the Cauchy problem for the nonlocal
fractional parabolic ADE:

\begin{equation}
\partial _{t}u+a\ast D_{x}^{\gamma }u+A\ast u=f\left( t,x\right) ,\text{ }%
u(0,x)=0,  \tag{1.4}
\end{equation}%
in $H$-valued mixed spaces $L_{\mathbf{p}}$ for $\mathbf{p=}\left(
p,p_{1}\right) $. In other words, we show that problem $\left( 1.4\right) $
has a unique solution $u\in W_{\mathbf{p}}^{1,\gamma }\left( \mathbb{R}%
_{+}^{2};H\left( A\right) ,E\right) $ for $f\in L_{\mathbf{p}}\left( \mathbb{%
R}_{+}^{2};H\right) $\ satisfying the following coercive estimate 
\begin{equation*}
\left\Vert \partial _{t}u\right\Vert _{L_{_{\mathbf{p}}}\left( \mathbb{R}%
_{+}^{2};H\right) }+\left\Vert D_{x}^{\gamma }u\right\Vert _{L_{\mathbf{p}%
}\left( \mathbb{R}_{+}^{2};H\right) }+\left\Vert A\ast u\right\Vert _{L_{%
\mathbf{p}};\left( \mathbb{R}_{+}^{2};H\right) }\leq 
\end{equation*}%
\begin{equation}
M\left\Vert f\right\Vert _{L_{\mathbf{p}}\left( \mathbb{R}_{+}^{2};H\right)
},  \tag{1.5}
\end{equation}%
here $L_{\mathbf{p}}=$ $L_{\mathbf{p}}\left( \mathbb{R}_{+}^{2};H\right) $
denote the space of $H-$valued strongly measurable functions $f$ \ defined\
on $\mathbb{R}_{+}^{2}$\ equipped with the mixed norm 
\begin{equation*}
\left\Vert f\right\Vert _{L_{\mathbf{p}}\left( \mathbb{R}_{+}^{2};H\right)
}=\left( \int\limits_{\mathbb{R}}\left( \int\limits_{0}^{\infty }\left\Vert
f\left( t,x\right) \right\Vert _{H}^{p_{1}}dt\right) ^{\frac{p}{p_{1}}%
}dx\right) ^{\frac{1}{p}}<\infty \text{, }p_{1},p\in \left( 1,\infty \right)
.
\end{equation*}%
As an application, in this paper the following are established: (a) maximal
regularity properties of the anisotropic elliptic fractional ADE in mixed $%
L_{\mathbf{p}}$, $\mathbf{p=}\left( p_{1},p\right) $ spaces; (b) coercive
properties of the system of FDEs of infinite many order in $L_{p}$ spaces.
Let $\mathbb{C}$ denote the set of complex numbers and 
\begin{equation*}
S_{\varphi }=\left\{ \lambda ;\text{ \ }\lambda \in \mathbb{C}\text{, }%
\left\vert \arg \lambda \right\vert \leq \varphi \right\} \cup \left\{
0\right\} ,\text{ }0\leq \varphi <\pi .\ 
\end{equation*}

$B\left( E_{1},E_{2}\right) $ denotes the space of bounded linear operators
from $E_{1}$ to $E_{2}$. For $E_{1}=E_{2}=E$ it denotes by $B\left( E\right) 
$. Let $D\left( A\right) $, $R\left( A\right) $ denote the domain and range
of the linear operator in $E,$ respectively. Let Ker $A$ denote a null space
of $A$. A closed linear operator\ $A$ is said to be $\varphi -$ sectorial
(or sectorial for $\varphi =0$) in a Banach\ space $E$ with bound $M>0$ if
Ker $A=\left\{ 0\right\} $, $D\left( A\right) $ and $R\left( A\right) $ are
dense on $E,$ and $\left\Vert \left( A+\lambda I\right) ^{-1}\right\Vert
_{B\left( E\right) }\leq M\left\vert \lambda \right\vert ^{-1}$ for all $%
\lambda \in S_{\varphi },$ $\varphi \in \left[ 0,\right. \left. \pi \right) $%
, where $I$ is an identity operator in $E.$ Sometimes $A+\lambda I$\ will be
written as $A+\lambda $ and will be denoted by $A_{\lambda }$. It is known $%
\left[ \text{20, \S 1.15.1}\right] $ that the powers\ $A^{\theta }$, $\theta
\in \left( -\infty ,\infty \right) $ for a positive operator $A$ exist$.$
Let $E\left( A^{\theta }\right) $ denote the space $D\left( A^{\theta
}\right) $ with the norm 
\begin{equation*}
\left\Vert u\right\Vert _{E\left( A^{\theta }\right) }=\left( \left\Vert
u\right\Vert ^{p}+\left\Vert A^{\theta }u\right\Vert ^{p}\right) ^{\frac{1}{p%
}},\text{ }1\leq p<\infty ,\text{ }0<\theta <\infty .
\end{equation*}

A sectorial operator $A\left( x\right) ,$ $x\in \mathbb{R}^{n}$ is said to
be uniformly sectorial in a Banach space $E$ if there exists a $\varphi \in %
\left[ 0\right. ,\left. \pi \right) $ such that the uniformly estimate holds 
\begin{equation*}
\left\Vert \left( A\left( x\right) +\lambda I\right) ^{-1}\right\Vert
_{B\left( E\right) }\leq M\left\vert \lambda \right\vert ^{-1}
\end{equation*}%
for all $\lambda \in S_{\varphi }.$

Let $S\left( \mathbb{R}^{n};E\right) $ denote the $E$-valued Schwartz class,
i.e., the space of all $E$-valued rapidly decreasing smooth functions on $%
\mathbb{R}^{n}$ equipped with its usual topology generated by seminorms. For 
$E=\mathbb{C}$ this space will be denoted by $S=S\left( \mathbb{R}%
^{n}\right) $. Here, $S^{\prime }\left( E\right) =S^{\prime }\left( \mathbb{R%
}^{n};E\right) $ denotes the space of linear continuous mappings from $S$
into\ $E$ and is called $E$-valued Schwartz distributions. For any $\alpha
=\left( \alpha _{1},\alpha _{2},...,\alpha _{n}\right) $, $\alpha _{i}\in %
\left[ 0,\infty \right) $ the function $\left( i\xi \right) ^{\alpha }$ will
be defined as:

{\large 
\begin{equation*}
\left( i\xi \right) ^{\alpha }=\left\{ 
\begin{array}{c}
\left( i\xi _{1}\right) ^{\alpha _{1}},.,\left( i\xi _{n}\right) ^{\alpha
_{n}}\text{, }\xi _{1}\xi _{2},.,\xi _{n}\neq 0 \\ 
0\text{, }\xi _{1},\xi _{2},.,\xi _{n}=0,%
\end{array}%
\right. 
\end{equation*}%
}where 
\begin{equation*}
\left( i\xi _{k}\right) ^{\alpha _{k}}=\exp \left[ \alpha _{k}\left( \ln
\left\vert \xi _{k}\right\vert +i\frac{\pi }{2}\text{ sgn }\xi _{k}\right) %
\right] \text{, }k=1,2,...,n.
\end{equation*}%
The Liouville derivatives $D^{\alpha }u$ of an $E$-valued function $u$ are
defined similarly to the case of scalar functions $\left[ \text{13}\right] $%
. 

$C\left( \Omega ;E\right) $ and $C^{\left( m\right) }\left( \Omega ;E\right) 
$\ will denote the spaces of $E$-valued bounded uniformly strongly
continuous and $m$ times continuously differentiable functions on $\Omega $,
respectively. Let $F$ and $F^{-1}$ denote the Fourier and inverse Fourier
transforms. Through this section, the Fourier transformation of a function $u
$ will be denoted by $\hat{u}$. Let $E_{1}$ and $E_{2}$ be two Banach
spaces. A function $\Psi \in C\left( \mathbb{R}^{n};B\left(
E_{1},E_{2}\right) \right) $ is called a Fourier multiplier from $%
L_{p}\left( \mathbb{R}^{n};E_{1}\right) $\ to $L_{p}\left( \mathbb{R}%
^{n};E_{2}\right) $ if the map%
\begin{equation*}
u\rightarrow \Lambda u=F^{-1}\Psi \left( \xi \right) Fu,\text{ }u\in S\left( 
\mathbb{R}^{n};E_{1}\right) 
\end{equation*}%
is well defined and extends to a bounded linear operator 
\begin{equation*}
\Lambda :\ L_{p}\left( \mathbb{R}^{n};E_{1}\right) \rightarrow \ L_{p}\left( 
\mathbb{R}^{n};E_{2}\right) .
\end{equation*}

Let $E_{0}$ and $E$ be two Banach spaces and $E_{0}$ be continuously and
densely embedded into $E$. Let $s\in \mathbb{R}$ and $\xi =\left( \xi
_{1},\xi _{2},...,\xi _{n}\right) \in \mathbb{R}^{n}$. Consider the
following Liouville-Lions space

{\large 
\begin{equation*}
W_{p}^{s}(\mathbb{R}^{n};E_{0},E)=\left\{ u\right. \ u\in S^{\prime }\left( 
\mathbb{R}^{n};E_{0}\right) \text{,  }
\end{equation*}%
}%
\begin{equation*}
F^{-1}\left( 1+\left\vert \xi \right\vert ^{2}\right) ^{\frac{s}{2}}Fu\in
L_{p}\left( \mathbb{R}^{n};E\right) \text{,}
\end{equation*}

\begin{equation*}
\left\Vert u\right\Vert _{W_{p}^{s}\left( \mathbb{R}^{n};E_{0},E\right)
}=\left\Vert u\right\Vert _{L_{p}\left( \mathbb{R}^{n};E_{0}\right) }+\text{ 
}\left. \left\Vert F^{-1}\left( 1+\left\vert \xi \right\vert ^{2}\right) ^{%
\frac{s}{2}}Fu\right\Vert _{L_{p}\left( \mathbb{R}^{n};E\right) }<\infty
\right\} .
\end{equation*}

Sometimes we use one and the same symbol $C$ without distinction in order to
denote positive constants which may differ from each other even in a single
context. When we want to specify the dependence of such a constant on a
parameter, say $\alpha $, we write $C_{\alpha }$.

The embedding theorems in vector valued spaces play a key role in the theory
of DOEs. By reasoning as in $\left[ 16\right] $ we obtin estimating lower
order derivatives in terms of interpolation spaces:

{\large \ }\textbf{Theorem A}$_{1}$\textbf{.} Suppose $H$ is a Hilbert
space, $1<p\leq q<\infty $ and $A$ is a positive operator in $H$. Then for $s
$ $\in \left( 0,\infty \right) $ with $\varkappa =\frac{1}{s}\left[
\left\vert \alpha \right\vert +n\left( \frac{1}{p}-\frac{1}{q}\right) \right]
\leq 1$, $0\leq \mu \leq 1-\varkappa $ the embedding 
\begin{equation*}
D^{\alpha }W_{p}^{s}\left( \mathbb{R}^{n};H\left( A\right) ,H\right) \subset
L_{q}\left( \mathbb{R}^{n};H\left( A^{1-\varkappa -\mu }\right) \right) 
\end{equation*}%
is continuous and there exists a constant \ $C_{\mu }$ \ $>0$, depending
only on $\mu $ such that 
\begin{equation*}
\left\Vert D^{\alpha }u\right\Vert _{L_{q}\left( \mathbb{R}^{n};H\left(
A^{1-\varkappa -\mu }\right) \right) }\leq C_{\mu }\left[ h^{\mu }\left\Vert
u\right\Vert _{W_{p}^{s}\left( \mathbb{R}^{n};H\left( A\right) ,H\right)
}+h^{-\left( 1-\mu \right) }\left\Vert u\right\Vert _{L_{p}\left( \mathbb{R}%
^{n};H\right) }\right] 
\end{equation*}%
for all $u\in W_{p}^{s}\left( \mathbb{R}^{n};H\left( A\right) ,H\right) $
and $0<h\leq h_{0}<\infty .$

\begin{center}
2. \textbf{Elliptic} \textbf{FADE with parameters }
\end{center}

\bigskip Consider the problem $\left( 1.1\right) $.

\textbf{Condition 2.1. }Assume $\hat{A}(\xi )$ is a uniformly $\varphi -$%
sectorial operator in $H$ for $\varphi \in \left[ 0,\right. \left. \pi
\right) $ and $\hat{a}\in C^{\left( 1\right) }\left( \mathbb{R}\right) $
such that 
\begin{equation}
\hat{a}\left( \xi \right) \left( i\xi \right) ^{\gamma }\in S_{\varphi _{1}}%
\text{ for all }\xi \in \mathbb{R}\text{ and }\left\vert \hat{a}\left( \xi
\right) \right\vert \left\vert \left( i\xi \right) ^{\gamma }\right\vert
\leq C_{0}\xi ^{2}  \tag{2.1}
\end{equation}

Suppose $\left[ D^{\beta }\hat{A}\left( \xi \right) \right] \hat{A}%
^{-1}\left( \xi _{0}\right) \in C\left( \mathbb{R};B\left( H\right) \right) $
and 
\begin{equation}
\text{ }\left\vert \xi \right\vert ^{\left\vert \beta \right\vert
}\left\vert D^{\beta }\hat{a}(\xi )\right\vert \leq C_{1}\text{, }\beta \in
\left\{ 0,1\right\} ,\text{ }\xi \in \mathbb{R}\backslash \left\{ 0\right\} ,
\tag{2.2}
\end{equation}%
\begin{equation}
\left\Vert \left\vert \xi \right\vert ^{\beta }\text{ }D^{\beta }\hat{A}%
\left( \xi \right) \hat{A}^{-1}\left( \xi _{0}\right) \right\Vert
_{B(H)}\leq C_{2}\text{, }\beta \in \left\{ 0,1\right\} ,\text{ }\xi ,\text{ 
}\xi _{0}\in \mathbb{R}\backslash \left\{ 0\right\} \text{.}  \tag{2.3}
\end{equation}

Let 
\begin{equation*}
X=L_{p}\left( \mathbb{R};H\right) \text{, }Y=W_{p}^{2}\left( \mathbb{R}%
;H\left( A\right) ,H\right) .
\end{equation*}%
\ In this section we prove the following:

\textbf{Theorem 2.1.} Assume that the Condition 2.1 is satisfied. Suppose
that $\gamma \in \left( \left. 1,2\right] ,\right. A$ is a sectorial
operator in $H$ for $\varphi \in \left( 0\right. ,\left. \pi \right] $ and $%
\lambda \in S_{\varphi _{2}}$. Then for $f\in X$, $0\leq \varphi _{1}<\pi
-\varphi _{2}$ and $\varphi _{1}+\varphi _{2}\leq \varphi $ there is a
unique solution $u$ of the equation $\left( 1.1\right) $ belonging to $Y$
and the following coercive uniform estimate holds 
\begin{equation}
\dsum\limits_{0\leq s\leq \gamma }\left\vert \lambda \right\vert ^{1-\frac{s%
}{\gamma }}\left\Vert a\ast D^{s}u\right\Vert _{X}+\left\Vert A\ast
u\right\Vert _{X}\leq C\left\Vert f\right\Vert _{X}.  \tag{2.4}
\end{equation}

\textbf{Proof. }By applying the Fourier transform to equation $\left(
1.1\right) $ we\ obtain 
\begin{equation}
Q\left( \xi ,\lambda \right) \hat{u}\left( \xi \right) =\hat{f}\left( \xi
\right) .  \tag{2.5}
\end{equation}%
where 
\begin{equation}
Q\left( \xi ,\lambda \right) =\hat{a}\left( \xi \right) \left[ \left( i\xi
\right) ^{\gamma }+\hat{A}\left( \xi \right) +\lambda \right] .  \tag{2.6}
\end{equation}%
By assumption, $\lambda +\hat{a}\left( \xi \right) \left( i\xi \right)
^{\gamma }\in S_{\varphi }$, for all $\xi \in \mathbb{R}$ and the operator\ $%
Q\left( \xi ,\lambda \right) $ is invertible in $H$ for all $\xi \in \mathbb{%
R}$ and $\lambda \in S_{\varphi _{2}}$. So, from $\left( 2.6\right) $ we
obtain that the solution of $\left( 2.5\right) $ can be represented in the
form 
\begin{equation}
u\left( x\right) =F^{-1}\left[ Q\left( \xi ,\lambda \right) \right] ^{-1}%
\hat{f}.  \tag{2.7}
\end{equation}%
By definition of the sectoriale operator $A$, the inverse of $A^{-1}$ is
bounded in $H$. Then the operator $A$ is a closed linear operator (as an
inverse of bounded linear operator $A^{-1}$). By properties of the Fourier
transform and by using $\left( 2.7\right) $ we have

{\large 
\begin{equation*}
\left\Vert Au\right\Vert _{X}=\left\Vert F^{-1}\left[ Q\left( \xi ,\lambda
\right) \right] ^{-1}\hat{f}\right\Vert _{X}\text{, }
\end{equation*}%
}%
\begin{equation*}
\left\Vert D^{s}u\right\Vert _{X}=\left\Vert F^{-1}\left[ \xi ^{s}Q\left(
\xi ,\lambda \right) \right] ^{-1}\hat{f}\right\Vert _{X}.
\end{equation*}%
Hence, it suffices to show that operator-functions 
\begin{equation*}
\eta \left( \lambda ,\xi \right) =A\left[ Q\left( \xi ,\lambda \right) %
\right] ^{-1}\text{, }\eta _{s}\left( \lambda ,\xi \right)
=\dsum\limits_{0\leq s\leq \gamma }\left\vert \lambda \right\vert ^{1-\frac{s%
}{\gamma }}\xi ^{s}\left[ Q\left( \xi ,\lambda \right) \right] ^{-1}
\end{equation*}%
are collections of multipliers in $X\ $uniformly with respect to $\lambda
\in S_{\varphi _{2}}$. By virtue of $\left[ \text{6, Lemma 2.3}\right] $,
for $\lambda \in S_{\varphi _{2}}$ and $\nu \in S_{\varphi _{1}}$ with $%
\varphi _{1}+\varphi _{2}<\pi $ there is a positive constant $C$ such that 
\begin{equation}
\left\vert \lambda +\nu \right\vert \geq C\left( \left\vert \lambda
\right\vert +\left\vert \nu \right\vert \right) .  \tag{2.8}
\end{equation}%
By using the resolvent properties of the operator $A$, we get that $\left[
Q\left( \xi ,\lambda \right) \right] ^{-1}$ is uniformly bounded for all $%
\xi \in \mathbb{R}$, $\lambda \in S_{\varphi _{2}}$ and 
\begin{equation*}
\left\Vert \left[ Q\left( \xi ,\lambda \right) \right] ^{-1}\right\Vert \leq
C\left( 1+\left\vert \lambda +\hat{a}\left( \xi \right) \left( i\xi \right)
^{\gamma }\right\vert \right) ^{-1}.
\end{equation*}%
By $\left( 2.1\right) $ and $\left( 2.8\right) $, we obtain that 
\begin{equation}
\left\Vert \left[ Q\left( \xi ,\lambda \right) \right] ^{-1}\right\Vert \leq
C\left( 1+\left\vert \lambda \right\vert +\left\vert \hat{a}\left( \xi
\right) \right\vert \left\vert \xi \right\vert ^{\gamma }\right) ^{-1}\leq  
\tag{2.9}
\end{equation}%
\begin{equation*}
C_{2}\left[ 1+\left\vert \lambda \right\vert +\xi ^{2}\right] ^{-1}.
\end{equation*}%
Then by resolvent properties of sectorial operators and uniform estimate $%
\left( 2.9\right) $\ we get 
\begin{equation*}
\left\Vert \eta \left( \lambda ,\xi \right) \right\Vert \leq \left\Vert
I+\left( \lambda +a\left( i\xi \right) ^{\gamma }\right) \left[ Q\left( \xi
,\lambda \right) \right] ^{-1}\right\Vert \leq 
\end{equation*}%
\begin{equation*}
1+\left( \left\vert \lambda \right\vert +\left\vert a\right\vert \left\vert
\xi \right\vert ^{\gamma }\right) \left( 1+\left\vert \lambda \right\vert
+\left\vert a\right\vert \left\vert \xi \right\vert ^{\gamma }\right)
^{-1}\leq C_{3},
\end{equation*}%
where $I$ is an identity operator in $E$. Moreover, by well-known
inequality, we have%
\begin{equation*}
\left\vert \lambda \right\vert ^{1-\frac{s}{2}}\left\vert \xi \right\vert
^{s}=\left\vert \lambda \right\vert \left\vert \lambda \right\vert ^{-\frac{s%
}{\gamma }}\left\vert \xi \right\vert ^{s}\leq \left\vert \lambda
\right\vert \left( \left\vert \lambda \right\vert ^{-\frac{1}{\gamma }%
}\left\vert \xi \right\vert \right) ^{s}\leq 
\end{equation*}%
\begin{equation}
\left\vert \lambda \right\vert \left[ 1+\left\vert \lambda \right\vert
^{-1}\left\vert \xi \right\vert ^{\gamma }\right] =\left\vert \lambda
\right\vert +\left\vert \xi \right\vert ^{\gamma }.  \tag{2.10}
\end{equation}

Hence, in view of $\left( 2.9\right) $ and $\left( 2.10\right) $, we have%
\begin{equation*}
\left\Vert \eta _{s}\left( \lambda ,\xi \right) u\right\Vert _{H}\leq
\dsum\limits_{0\leq s\leq \gamma }\left\vert \lambda \right\vert ^{1-\frac{s%
}{\gamma }}\xi ^{s}\left\Vert \left[ Q\left( \xi ,\lambda \right) \right]
^{-1}u\right\Vert _{H}\leq C_{2}\left\Vert u\right\Vert _{H}.
\end{equation*}%
\ By using $\left( 2.2\right) $ and $\left( 2.3\right) $\ we obtain that the
operator functions $D^{i}\eta \left( \lambda ,\xi \right) $ and $D^{i}\eta
_{s}\left( \lambda ,\xi \right) $ are uniformly bounded, i.e., 
\begin{equation*}
\left\Vert D^{i}\eta \left( \lambda ,\xi \right) \right\Vert _{B\left(
H\right) }\leq C_{1},\left\Vert D^{i}\eta _{s}\left( \lambda ,\xi \right)
\right\Vert _{B\left( E\right) }\leq C_{2}
\end{equation*}%
for $i=0,1.$\ Then by virtue of $\left[ \text{21, Theorem 4.1}\right] $, $%
\eta \left( \lambda ,\xi \right) $ and $\eta _{s}\left( \lambda ,\xi \right) 
$ are Fourier multipliers in $L_{p}\left( \mathbb{R};H\right) .$

Let $O$ denote the operator in $X$ generated by problem $\left( 1.1\right) $
for $\lambda =0$, i.e., 
\begin{equation*}
D\left( O\right) \subset W_{p}^{\gamma }\left( \mathbb{R}^{n};H\left(
A\right) ,H\right) ,\text{ }Ou=a\ast D^{\gamma }u+A\ast u.
\end{equation*}

Theorem 2.1 and the definition of the space $W_{p}^{s}\left( \mathbb{R}%
^{n};H\left( A\right) ,H\right) $ imply the following result:

\textbf{Result 2.1.} Theorem 2.1 implies that the operator $O$ is separable
in $X$, i.e. for all $f\in X$ \ there is a unique solution $u\in Y$\ of the
problem $\left( 1.1\right) $, all terms of equation $\left( 1.1\right) $ are
also from $X$ and there are positive constants $C_{1}$ and $C_{2}$ so that 
\begin{equation}
C_{1}\left\Vert Ou\right\Vert _{X}\leq \dsum\limits_{0\leq s\leq \gamma
}\left\Vert a\ast D^{s}u\right\Vert _{X}+\left\Vert A\ast u\right\Vert
_{X}\leq C_{2}\left\Vert Ou\right\Vert _{X}.  \tag{2.11}
\end{equation}

Indeed, if we put $\lambda =1$ in $\left( 2.4\right) ,$ by Theorem 2.1 we
get the second inequality. So it is remain to prove the first estimate.\ The
first inequality is equivalent to the following estimate%
\begin{equation*}
\left\Vert F^{-1}\hat{A}\hat{u}\right\Vert _{X}+\left\Vert F^{-1}\hat{a}%
\left( i\xi \right) ^{\gamma }\hat{u}\right\Vert _{X}\leq
\end{equation*}

\begin{equation*}
C\left\{ \left\Vert F^{-1}\hat{A}\hat{u}\right\Vert _{X}+\left\Vert F^{-1}%
\left[ 1+\left( \dsum\limits_{k=1}^{n}\xi _{k}^{2}\right) ^{\frac{1}{2}}%
\right] ^{2}\hat{u}\right\Vert _{X}\right\} .
\end{equation*}%
So, it suffices to show that the operator functions%
\begin{equation*}
A\left\{ A+\left[ 1+\left( \dsum\limits_{k=1}^{n}\xi _{k}^{2}\right) ^{\frac{%
1}{2}}\right] ^{2}\right\} ^{-1},\text{ }\hat{a}\left( \xi \right) \xi
^{\gamma }\left[ 1+\left( \dsum\limits_{k=1}^{n}\xi _{k}^{2}\right) ^{\frac{1%
}{2}}\right] 
\end{equation*}%
are uniform Fourier multipliers in $X$. This fact is proved in a similar way
as in the proof of Theorem 2.1.

From Theorem 2.1, we have:

\textbf{Result 2.2. }Assume all conditions of Theorem 2.1 hold. Then, for
all $\lambda \in S_{\varphi }$ the resolvent of operator $O$ exists and the
following sharp coercive uniform estimate holds 
\begin{equation}
\dsum\limits_{0\leq s\leq \gamma }\left\vert \lambda \right\vert ^{1-\frac{s%
}{\gamma }}\left\Vert a\ast D^{s}\left( O+\lambda \right) ^{-1}\right\Vert
_{B\left( X\right) }+\left\Vert A\ast \left( O+\lambda \right)
^{-1}\right\Vert _{B\left( X\right) }\leq C.  \tag{2.11}
\end{equation}

Indeed, we infer from Theorem 2.1 that the operator $O+\lambda $ has a
bounded inverse from $X$ to $Y.$ So, the solution\ $u$\ of the equation $%
\left( 1.1\right) $ can be expressed as $u\left( x\right) =\left( O+\lambda
\right) ^{-1}f$ \ for all $f\in X.$ Then estimate $\left( 2.4\right) $
implies the estimate $\left( 2.11\right) .$

\textbf{Condition 2.2.} Let $D(A(x))=D(\hat{A}(\xi ))$, $D(\hat{A}(\xi ))$
is dense in $E$ and does not depend on $\xi ;$ $A\left( x\right) $ is a
uniformly sectorial in $H.$ Moreover, there exist positive constants $C_{1}$%
, $C_{2}$ and $\xi _{0}\in \mathbb{R}$ such that

\begin{equation*}
C_{1}\left\Vert \hat{A}\left( \xi _{0}\right) u\right\Vert _{E}\leq
\left\Vert A\left( x\right) u\right\Vert _{E}\leq C_{2}\left\Vert \hat{A}%
\left( \xi _{0}\right) u\right\Vert _{E}
\end{equation*}%
for $u\in D\left( A\right) $, $x\in \mathbb{R}.$

Then, ee prove the following

\textbf{Theorem 2.2.} Assume that the Conditions 2.1, 2.2 are satisfied.
Suppose that $\gamma \in \left( \left. 1,2\right] ,\right. A$ is a sectorial
operator in $H$ with respect to $\varphi \in \left( 0\right. ,\left. \pi %
\right] $ and $\lambda \in S_{\varphi _{2}}$. Then for $f\in X$, $0\leq
\varphi _{1}<\pi -\varphi _{2}$ and $\varphi _{1}+\varphi _{2}\leq \varphi $
there is a unique solution $u$ of the equation $\left( 1.1\right) $
belonging to $Y$ and the following coercive uniform estimate holds 
\begin{equation}
\dsum\limits_{0\leq s\leq \gamma }\left\vert \lambda \right\vert ^{1-\frac{s%
}{\gamma }}\left\Vert D^{s}u\right\Vert _{X}+\left\Vert Au\right\Vert
_{X}\leq C\left\Vert f\right\Vert _{X}.  \tag{2.12}
\end{equation}

\textbf{Proof}. The estimate $\left( 2.12\right) $ is derived by reasoning
as in Theorem 2.2.

From Theorem 2.2, we have the following results:

\textbf{Result 2.3.} There are positive constants $C_{1}$ and $C_{2}$ so
that 
\begin{equation}
C_{1}\left\Vert Ou\right\Vert _{X}\leq \dsum\limits_{0\leq s\leq \gamma
}\left\Vert D^{s}u\right\Vert _{X}+\left\Vert Au\right\Vert _{X}\leq
C_{2}\left\Vert Ou\right\Vert _{X}.  \tag{2.13}
\end{equation}

Indeed, the second inequality $\left( 2.13\right) $ is obtained from Theorem
2.1. So, it remains to prove the first estimate. Indeed, from Condition 2.2
we have 
\begin{equation*}
\left\Vert A\ast u\right\Vert _{X}\leq M\left\Vert F^{-1}\hat{A}\hat{u}%
\right\Vert _{X}\leq C\left\Vert F^{-1}\hat{A}A^{-1}\left( x_{0}\right)
A\left( x_{0}\right) \hat{u}\right\Vert _{X}
\end{equation*}%
\begin{equation}
\leq C\left\Vert F^{-1}A\left( x_{0}\right) \hat{u}\right\Vert _{X}\leq
C\left\Vert Au\right\Vert _{X}\text{, for }u\in Y.  \tag{2.14}
\end{equation}

Then by using $\left( 2.14\right) $ and by reasoning as in Theorem 2.1 we
obtain $\left( 2.13\right) .$

\textbf{Result 2.4. }Assume all conditions of Theorem 2.2 hold. Then, for
all $\lambda \in S_{\varphi }$ the resolvent of operator $O$ exists and the
following sharp uniform estimate holds 
\begin{equation}
\dsum\limits_{0\leq s\leq \gamma }\left\vert \lambda \right\vert ^{1-\frac{s%
}{\gamma }}\left\Vert D^{s}\left( O+\lambda \right) ^{-1}\right\Vert
_{B\left( X\right) }+\left\Vert \left( O+\lambda \right) ^{-1}\right\Vert
_{B\left( X\right) }\leq C.  \tag{2.15}
\end{equation}

\textbf{Result 2.5. }Theorem 2.2 particularly implies that the operator $O$
is sectorial in $X.$ Then the operators $O^{s}$ are generators of analytic
semigroups in $X$ for $s\leq \frac{1}{2}$ (see e.g. $\left[ \text{19, \S %
1.14.5}\right] $)$.$

\begin{center}
\bigskip \textbf{3. The Cauchy problem for fractional parabolic ADE }
\end{center}

In this section, we shall consider the following Cauchy problem for the
parabolic FDOE 
\begin{equation}
\frac{\partial u}{\partial t}+a\ast D_{x}^{\gamma }u+A\ast u=f\left(
t,x\right) ,\text{ }u(0,x)=0\text{, }t\in \mathbb{R}_{+}\text{, }x\in 
\mathbb{R},  \tag{3.1}
\end{equation}%
where $a$ is a complex number, $D_{x}^{\alpha }$ is the fractional
derivative in $x$ for $\alpha \in \left( 1,2\right) $ defined by $\left(
1.2\right) $ and $A$ is a linear operator in $E$.

By applying Theorem 2.1 we establish the maximal regularity of the problem $%
\left( 3.1\right) $ in $E$-valued mixed $L_{\mathbf{p}}$ spaces, where $%
\mathbf{p=}\left( p_{1},p\right) $. Let $O$ denote the operator generated by
problem $\left( 1.1\right) $ for $\lambda =0$. Let $E$ be a Banach space.
For $\mathbf{p=}\left( p\text{, }p_{1}\right) ,$ $Z=L_{\mathbf{p}}\left( 
\mathbb{R}_{+}^{2};H\right) $ will denote the space of all $\mathbf{p}$%
-summable $E$-valued\ functions on $\mathbb{R}_{+}^{2}$ with mixed norm,
i.e., the space of all measurable $E$-valued functions $f$ defined on for
which 
\begin{equation*}
\left\Vert f\right\Vert _{L_{\mathbf{p}}\left( \mathbb{R}_{+}^{2};H\right)
}=\left( \int\limits_{\mathbb{R}}\left( \int\limits_{0}^{\infty }\left\Vert
f\left( t,x\right) \right\Vert _{H}^{p}dx\right) ^{\frac{p_{1}}{p}}dt\right)
^{\frac{1}{p_{1}}}<\infty .
\end{equation*}

Let  $Z^{1,\gamma }\left( A\right) =W_{\mathbf{p}}^{1,\gamma }\left( \mathbb{%
R}_{+}^{2};H\left( A\right) ,H\right) $ denotes the space of all functions $%
u\in L_{\mathbf{p}}\left( \mathbb{R}_{+}^{n+1};H\left( A\right) \right) $
possessing the generalized derivative $D_{t}u=\frac{\partial u}{\partial t}%
\in Z$ with respect to $y$\ and fractional derivatives $D_{x}^{\gamma }u\in Z
$ with respect to $x$ for $\left\vert \alpha \right\vert \leq m$ with the
norm 
\begin{equation*}
\ \left\Vert u\right\Vert _{Z^{1,2}\left( A\right) }=\left\Vert
Au\right\Vert _{Z}+\left\Vert \partial _{t}u\right\Vert _{Z}+\left\Vert
D_{x}^{\gamma }u\right\Vert _{Z},
\end{equation*}%
where $\gamma \in \left( \left. 1,2\right] \right. $ and $u=u\left(
t,x\right) .$

Now, we are ready to state the main result of this section.

\textbf{Theorem 3.1.}\ Assume the conditions of Theorem 2.1 hold for $%
\varphi \in \left( \frac{\pi }{2},\pi \right) $. Then for $f\in Z$ problem $%
\left( 3.1\right) $ has a unique solution $u\in Z^{1,\gamma }\left( A\right) 
$ satisfying the following unform coercive estimate 
\begin{equation*}
\left\Vert \partial _{t}u\right\Vert _{Z}+\left\Vert a\ast D_{x}^{\gamma
}u\right\Vert _{Z}+\left\Vert A\ast u\right\Vert _{Z}\leq C\left\Vert
f\right\Vert _{Z}.
\end{equation*}

\textbf{Proof.} By definition of $X=L_{p}\left( \mathbb{R};H\right) $ and
mixed space $L_{\mathbf{p}}\left( \mathbb{R}_{+}^{2};H\right) ,$ $\mathbf{p=}%
\left( p\text{, }p_{1}\right) $, we have

\begin{equation*}
\left\Vert u\right\Vert _{L_{p_{1}}\left( 0,\infty ;X\right) }=\left(
\dint\limits_{0}^{\infty }\left\Vert u\left( t\right) \right\Vert
_{X}^{p_{1}}dt\right) ^{\frac{1}{p_{1}}}=\left( \dint\limits_{0}^{\infty
}\left\Vert u\left( t\right) \right\Vert _{L_{p}\left( \mathbb{R};E\right)
}^{p_{1}}dt\right) ^{\frac{1}{p_{1}}}=\left\Vert u\right\Vert _{Z}.
\end{equation*}

Therefore, the problem $\left( 3.1\right) $\ can be expressed as the
following Cauchy problem for the abstract parabolic equation 
\begin{equation}
\frac{du}{dt}+Ou\left( t\right) =f\left( t\right) ,\text{ }u\left( 0\right)
=0,\text{ }t\in \left( 0,\infty \right) .  \tag{3.2}
\end{equation}

Then, by virtue of [21, Theorem 4.2], we obtain that for $f\in
L_{p_{1}}\left( 0,\infty ;X\right) $ the problem $\left( 3.2\right) $ has a
unique solution $u\in W_{p_{1}}^{1}\left( 0,\infty ;D\left( O\right)
,X\right) $ satisfying the following estimate 
\begin{equation*}
\left\Vert \frac{du}{dt}\right\Vert _{L_{p_{1}}\left( 0,\infty ;X\right)
}+\left\Vert Ou\right\Vert _{L_{p_{1}}\left( 0,\infty ;X\right) }\leq
C\left\Vert f\right\Vert _{L_{p_{1}}\left( 0,\infty ;X\right) }.
\end{equation*}

From the Theorem 2.2, relation $\left( 3.2\right) $ and from the above
estimate we get the assertion$.$

\begin{center}
\textbf{4. BVP for anisotropic elliptic FDE }
\end{center}

In this section, the maximal regularity properties of the anisotropic FDE is
studied. Let $\tilde{\Omega}=\mathbb{R\times }\Omega $, where $\Omega
\subset \mathbb{R}^{n}$ is an open connected set with compact $C^{2l}$%
-boundary $\partial \Omega $. Consider the BVP for the FDE

\begin{equation}
a\ast D_{x}^{\gamma }u+\sum\limits_{\left\vert \beta \right\vert \leq
2l}b_{\beta }\ast D_{y}^{\beta }+\lambda u=f\left( x,y\right) ,\text{ }y\in
\Omega ,  \tag{4.1}
\end{equation}%
\ \ \ 

\begin{equation}
B_{j}u=\sum\limits_{\left\vert \beta \right\vert \leq l_{j}}\ b_{j\beta
}\left( y\right) D_{y}^{\beta }u\left( x,y\right) =0\text{, }x\in \mathbb{R}%
\text{,}  \tag{4.2}
\end{equation}

\begin{equation*}
\text{ }y\in \partial \Omega ,\text{ }j=1,2,...,l,
\end{equation*}%
where $u=u\left( x,y\right) $, $a$ is a complex number,$\ $\ $D_{x}^{\gamma }
$ is the fractional derivative operator in $x$ for $\gamma \in \left( \left.
1,2\right] \right. $ defined by $\left( 1.2\right) $, 
\begin{equation*}
D_{j}=-i\frac{\partial }{\partial y_{j}}\text{, }y=\left(
y_{1},...,y_{n}\right) \text{, }b_{\alpha }=b_{\alpha }\left( y\right) \text{%
,}
\end{equation*}%
$\beta =\left( \beta _{1},\beta _{2},...,\beta _{n}\right) $ are nonnegative
integer numbers and $\lambda $ is a complex parameter. For $\mathbf{p=}%
\left( 2,p\right) $ here, $L_{\mathbf{p}}\left( \tilde{\Omega}\right) $ will
denote the space of all $\mathbf{p}$-summable scalar-valued\ functions with
mixed norm i.e., the space of all measurable functions $f$ defined on $%
\tilde{\Omega}$, for which 
\begin{equation*}
\left\Vert f\right\Vert _{L_{\mathbf{p}}\left( \tilde{\Omega}\right)
}=\left( \int\limits_{\mathbb{R}}\left( \int\limits_{\Omega }\left\vert
f\left( x,y\right) \right\vert ^{2}dx\right) ^{\frac{2}{p_{1}}}dy\right) ^{%
\frac{1}{2}}<\infty .
\end{equation*}%
Analogously, $W_{\mathbf{p}}^{\gamma ,2l}\left( \tilde{\Omega}\right) $
denotes the anisotropic fractional Sobolev space with corresponding mixed
norm, i.e., $W_{\mathbf{p}}^{\gamma ,2l}\left( \tilde{\Omega}\right) $
denotes the space of all functions $u\in L_{\mathbf{p}}\left( \tilde{\Omega}%
\right) $ possessing the fractional derivatives $D_{x}^{\gamma }u\in L_{%
\mathbf{p}}\left( \tilde{\Omega}\right) $ with respect to $x$ for $\alpha
\in \left( 1,2\right) $ and generalized derivative $\frac{\partial ^{2l}u}{%
\partial y_{k}^{2l}}\in L_{\mathbf{p}}\left( \tilde{\Omega}\right) $ with
respect to $y$ with the norm 
\begin{equation*}
\ \left\Vert u\right\Vert _{W_{\mathbf{p}}^{\gamma ,2l}\left( \tilde{\Omega}%
\right) }=\left\Vert u\right\Vert _{L_{\mathbf{p}}\left( \tilde{\Omega}%
\right) }+\left\Vert a\ast D_{x}^{\gamma }u\right\Vert _{L_{\mathbf{p}%
}\left( \tilde{\Omega}\right) }+\dsum\limits_{k=1}^{n}\left\Vert \frac{%
\partial ^{2l}u}{\partial y_{k}^{2l}}\right\Vert _{L_{\mathbf{p}}\left( 
\tilde{\Omega}\right) }.
\end{equation*}

\ Let $Q$ denote the operator generated by problem $\left( 4.1\right)
-\left( 4.2\right) $, i.e., 
\begin{equation*}
D\left( Q\right) =W_{\mathbf{p}}^{\gamma ,2l}\left( \tilde{\Omega}%
,B_{j}\right) =\left\{ u:u\in W_{\mathbf{p}}^{\gamma ,2l}\left( \tilde{\Omega%
}\right) ,\text{ }B_{j}u=0,\text{ }j=1,2,...l\right\} ,
\end{equation*}%
\begin{equation*}
Qu=aD_{x}^{\gamma }u+bD_{x}^{\nu }u+\sum\limits_{\left\vert \alpha
\right\vert \leq 2l}b_{\alpha }D_{y}^{\alpha }u.
\end{equation*}

Let $\xi ^{\prime }=\left( \xi _{1},\xi _{2},...,\xi _{n-1}\right) \in 
\mathbb{R}^{\mu -1},$ $\beta ^{\prime }=\left( \beta _{1},\beta
_{2},...,\beta _{n-1}\right) \in Z^{\mu }$ and 
\begin{equation*}
\text{ }A\left( y_{0},\xi ^{\prime },D_{y}\right) =\sum\limits_{\left\vert
\beta ^{\prime }\right\vert +j\leq 2l}a_{\beta ^{\prime }}\left(
y_{0}\right) \xi _{1}^{\beta _{1}}\xi _{2}^{\beta _{2}},.,\xi _{n-1}^{\beta
_{n-1}}D_{n}^{j}\text{ for }y_{0}\in \bar{G}
\end{equation*}%
\begin{equation*}
B_{j}\left( y_{0},\xi ^{\prime },D_{y}\right) =\sum\limits_{\left\vert \beta
^{\prime }\right\vert +j\leq l_{j}}b_{j\beta ^{\prime }}\left( y_{0}\right)
\xi _{1}^{\beta _{1}}\xi _{2}^{\beta _{2}},.,\xi _{\mu -1}^{\beta _{\mu
-1}}D_{\mu }^{j}\text{ for }y_{0}\in \partial G.
\end{equation*}

\textbf{Condition 4.1. }Let the following conditions be satisfied:

(0) the Condition 2.1. holds;

(1) $b_{\beta }\in C\left( \bar{\Omega}\right) $ for each $\left\vert \beta
\right\vert =2l$ and $b_{\alpha }\in L_{\infty }\left( \Omega \right)
+L_{r_{k}}\left( \Omega \right) $ for each $\left\vert \alpha \right\vert
=k<2l$ with $r_{k}\geq p_{1}$, $p_{1}\in \left( 1,\infty \right) $ and $2l-k>%
\frac{l}{r_{k}};$

(2) $b_{j\beta }\in C^{2l-l_{j}}\left( \partial \Omega \right) $ for each $j$%
, $\beta $, $l_{j}<2l$, $p\in \left( 1,\infty \right) ,$ $\lambda \in
S_{\varphi }$, $\varphi \in \lbrack 0,\pi );$

(3) for $y\in \bar{\Omega}$, $\xi \in R^{n}$, $\sigma \in S_{\varphi _{0}}$, 
$\varphi _{0}\in \left( 0,\frac{\pi }{2}\right) $, $\left\vert \xi
\right\vert +\left\vert \sigma \right\vert \neq 0$ let $\sigma
+\sum\limits_{\left\vert \alpha \right\vert =2l}b_{\alpha }\left( y\right)
\xi ^{\alpha }\neq 0;$

(4) for each $y_{0}\in \partial \Omega $ local BVP in local coordinates
corresponding to $y_{0},$ 
\begin{equation*}
\lambda +A\left( y_{0},\xi ^{\prime },D_{y}\right) \vartheta \left( y\right)
=0,
\end{equation*}

\begin{equation*}
B_{j}\left( y_{0},\xi ^{\prime },D_{y}\right) \vartheta \left( 0\right)
=h_{j}\text{, }j=1,2,...,l
\end{equation*}%
has a unique solution $\vartheta \in C_{0}\left( 0,\infty \right) $ for all $%
h=\left( h_{1},h_{2},...,h_{l}\right) \in \mathbb{C}^{l}$ and for $\xi
^{^{\prime }}\in \mathbb{R}^{n-1}$.

\textbf{\ }Suppose $\nu =\left( \nu _{1},\nu _{2},...,\nu _{n}\right) $ are
nonnegative real numbers. In this section, we present the following result:

\textbf{Theorem 4.1}. Assume $\gamma \in \left( \left. 1,2\right] \text{ }%
\right. $\ and $\lambda \in S_{\varphi _{2}}$ for $0\leq \varphi _{1}<\pi
-\varphi _{2}$. Moreover, suppose Condition 4.1 are satisfied. Then for $\
f\in L_{\mathbf{p}}\left( \tilde{\Omega}\right) $, $\lambda \in S_{\varphi
}, $ $\varphi \in \left( 0,\right. \left. \pi \right] $, $\varphi
_{1}+\varphi _{2}\leq \varphi $ problem $\left( 4.1\right) -\left(
4.2\right) $ has a unique solution $u\in W_{p}^{2,2l}\left( \tilde{\Omega}%
\right) $ and the following coercive uniform estimate holds 
\begin{equation*}
\dsum\limits_{0\leq s\leq \gamma }\left\vert \lambda \right\vert ^{1-\frac{s%
}{\gamma }}\left\Vert a\ast D_{x}^{s}u\right\Vert _{L_{\mathbf{p}}\left( 
\tilde{\Omega}\right) }+\sum\limits_{\left\vert \beta \right\vert \leq
2l}\left\Vert D_{y}^{\beta }u\right\Vert _{L_{\mathbf{p}}\left( \tilde{\Omega%
}\right) }\leq C\left\Vert f\right\Vert _{L_{\mathbf{p}}\left( \tilde{\Omega}%
\right) }.
\end{equation*}

\ \textbf{Proof.} Let $H=L_{2}\left( \Omega \right) $. It is known $\left[ 4%
\right] $\ that $L_{p_{1}}\left( \Omega \right) $ is an $UMD$ space for $%
p_{1}\in \left( 1,\infty \right) $. Consider the operator $A$ defined by 
\begin{equation*}
D\left( A\right) =W_{p_{1}}^{2l}\left( \Omega ;B_{j}u=0\right) ,\text{ }%
Au=\sum\limits_{\left\vert \beta \right\vert \leq 2l}b_{\beta }\left(
x\right) D^{\beta }u\left( y\right) .
\end{equation*}

Therefore, the problem $\left( 4.1\right) -\left( 4.2\right) $ can be
rewritten in the form of $\left( 2.1\right) $, where $u\left( x\right)
=u\left( x,.\right) ,$ $f\left( x\right) =f\left( x,.\right) $\ are
functions with values in $H=L_{2}\left( \Omega \right) $. From $\left[ \text{%
5, Theorem 8.2}\right] $ we get that the following problem

\begin{equation}
\eta u\left( y\right) +\sum\limits_{\left\vert \beta \right\vert \leq
2l}b_{\beta }\left( y\right) D^{\beta }u\left( y\right) =f\left( y\right) 
\text{, }  \tag{4.3}
\end{equation}%
\begin{equation*}
B_{j}u=\sum\limits_{\left\vert \beta \right\vert \leq l_{j}}\ b_{j\beta
}\left( y\right) D^{\beta }u\left( y\right) =0\text{, }j=1,2,...,l
\end{equation*}%
has a unique solution for $f\in L_{2}\left( \Omega \right) $ and arg $\eta
\in S\left( \varphi _{1}\right) ,$ $\left\vert \eta \right\vert \rightarrow
\infty .$ Moreover, the operator $A$ generated\ by $\left( 4.3\right) $ is
positive in $L_{2}$. Then from Theorem 2.1 we obtain the assertion.

\ \ \ \ \ 

\begin{center}
\textbf{5. The system of parabolic FDE of infinite many order }
\end{center}

Consider the following system of FDEs

\begin{equation}
\partial _{t}u_{i}+a\ast \partial _{x}^{\gamma
}u_{i}+\sum\limits_{j=1}^{N}\left( a_{ij}+\lambda \right) u_{j}\left(
x\right) =f_{i}\left( t,x\right) \text{, }t\text{, }x\in \mathbb{R}_{+}\text{%
, }  \tag{5.1}
\end{equation}

\begin{equation*}
u_{i}\left( 0,x\right) =0\text{, }i=1,2,...,N,\text{ }N\in \left[ 1,\infty %
\right] ,
\end{equation*}%
where $a$ is a complex number, $a_{ij}$ are real numbers and $\partial
_{x}^{\gamma }$ is the fractional differential operator in $x$ \ for $\gamma
\in \left( \left. 1,2\right] \text{ }\right. $ defined by $\left( 1.2\right) 
$. Let 
\begin{equation*}
\text{ }l_{2}\left( A\right) =\left\{ u\in l_{2},\left\Vert u\right\Vert
_{l_{q}\left( A\right) }=\left\Vert Au\right\Vert _{l_{2}}=\right.
\end{equation*}

\begin{equation*}
\left( \sum\limits_{i=1}^{N}\left\vert \left( Au\right) _{i}\right\vert
^{2}\right) ^{\frac{1}{2}}=\left. \left( \sum\limits_{i=1}^{N}\left\vert
\sum\limits_{j=1}^{N}a_{ij}u_{j}\right\vert ^{2}\right) ^{\frac{1}{2}%
}<\infty \right\} ,
\end{equation*}%
\begin{equation*}
\text{ }u=\left\{ u_{j}\right\} ,\text{ }Au=\left\{
\dsum\limits_{j=1}^{N}a_{ij}u_{j}\right\} ,\text{ }i,\text{ }j=1,2,...N.
\end{equation*}

\textbf{Condition 5.1. }Let the Condition 2.1. holds and 
\begin{equation*}
a_{ij}=a_{ji}\text{, }\sum\limits_{i,j=1}^{N}a_{ij}\xi _{i}\xi _{j}\geq
C_{0}\left\vert \xi \right\vert ^{2}\text{ for }\xi \neq 0.
\end{equation*}

Let%
\begin{equation*}
f\left( x\right) =\left\{ f_{i}\left( x\right) \right\} _{1}^{N}\text{, }%
u=\left\{ u_{i}\left( x\right) \right\} _{1}^{N},\text{ }Y_{2}=L_{\mathbf{p}%
}\left( \mathbb{R}_{+}^{2};l_{2}\right)
\end{equation*}

\textbf{Theorem 5.1. }Suppose the Condition 5.1 is satisfied.\textbf{\ }%
Then, for $f\left( t,x\right) \in L_{\mathbf{p}}\left( \mathbb{R}%
_{+}^{2};l_{2}\right) $ problem $\left( 5.1\right) $ has a unique solution $%
u\in $ $W_{\mathbf{p}}^{\gamma ,2}\left( \mathbb{R}_{+}^{2},l_{2}\left(
A\right) ,l_{2}\right) $ and the following uniform coercive estimate holds 
\begin{equation*}
\left\Vert \partial _{t}u\right\Vert _{Y_{2}}+\left\Vert a\ast D_{x}^{\gamma
}u\right\Vert _{Y_{2}}+\left\Vert u\right\Vert _{Y_{2}}\leq C\left\Vert
f\right\Vert _{Y_{2}}
\end{equation*}

\ \textbf{Proof.} Let $H=l_{2}$ and $A$ be a matrix such that $A=\left[
a_{ij}\right] $, $i$, $j=1,2,...,N.$ It is easy to see that 
\begin{equation*}
B\left( \lambda \right) =\lambda \left( A+\lambda \right) ^{-1}=\frac{%
\lambda }{D\left( \lambda \right) }\left[ A_{ji}\left( \lambda \right) %
\right] \text{, }i\text{, }j=1,2,...N,
\end{equation*}%
where $D\left( \lambda \right) =\det \left( A-\lambda I\right) $, $%
A_{ji}\left( \lambda \right) $ are entries of the corresponding adjoint
matrix of $A-\lambda I.$ Since the matrix $A$ is symmetric and positive
definite, it generates a positive operator in $l_{2}.$i.e., the operator $A$
is positive in $l_{q}$. So, from Theorem 3.1, we obtain the asserton.

From the above estimate we obtain the assertion.

\vspace{3mm}

\vspace{3mm}

\textbf{References}

\ \ \ \ \ \ \ \ \ \ \ \ \ \ \ \ \ \ \ \ \ \ \ \ \ \ \ \ \ \ \ \ \ \ \ \ \ \
\ \ \ \ \ \ \ \ \ \ \ \ \ \ \ \ \ \ \ \ \ \ \ \ \ \ \ \ \ \ \ \ \ \ \ \ \ \
\ \ \ \ \ \ \ \ \ 

\begin{enumerate}
\item H. Amann, Operator-valued Fourier multipliers, vector-valued Besov
spaces, and applications, Math. Nachr. 186 (1997), 5-56.

\item Agarwal R., Bohner M., Shakhmurov V. B., Maximal regular boundary
value problems in Banach-valued weighted spaces, Bound. Value Probl.,
(1)(2005), 9-42

\item Ashyralyev A, SharifovY.A., Existence and uniqueness of solutions for
the system of nonlinear fractional differential equations with nonlocal and
integral boundary conditions, Abstr. Appl. Anal., ID 594802, (2012).

\item Baleanu, D, Mousalou, A, Rezapour, Sh: A new method for investigating
some fractional integro-differential equations involving the Caputo-Fabrizio
derivative. Adv. Differ. Equ. 2017, 51 (2017).

\item Denk R., Hieber M., Pr\"{u}ss J., $R$-boundedness, Fourier multipliers
and problems of elliptic and parabolic type, Mem. Amer. Math. Soc. 166
(2003), n.788.

\item Dore C. and Yakubov S., Semigroup estimates and non coercive boundary
value problems, Semigroup Forum 60 (2000), 93-121.

\item Kilbas, A.A., Srivastava, H.M., Trujillo, J.J.: Theory and
Applications of Fractional Differential Equations. North-Holland Mathematics
Studies, vol. 204. Elsevier, Amsterdam (2006).

\item Favini A., Shakhmurov V. B, Yakov Yakubov, Regular boundary value
problems for complete second order elliptic differential-operator equations
in UMD Banach spaces, Semigroup form, 79 (1)( 2009), 22-54.

\item Krein S. G., Linear differential equations in Banach space, American
Mathematical Society, Providence, 1971.

\item Miller K. S. and Ross B., An Introduction to the Fractional Calculus
and Fractional Differential Equations, JohnWiley \& Sons, New York, NY, USA,
1993.

\item Lakshmikantham V., Vatsala A. S., Basic theory of fractional
differential equations, Nonlinear Analysis: Theory, Methods \& Applications,
69(8, (2008), 2677-2682.

\item Lakshmikantham V., Vasundhara J. D., Theory of fractional differential
Equations in a Banach space, European journal pure and applied mathematics,
(1)1, (2008), 38-45.

\item Lizorkin P. I., Generalized Liouville differentiation and functional
spaces $L_{p}^{r}\left( E_{n}\right) .$ Embedding theorems, Mathematics of
the USSR-Sbornik (3)60 (1963), 325-353.

\item Podlubny, I. Fractional Differential Equations San Diego: Academic
Press, 1999.

\item V. B. Shakhmurov, R. V. Shahmurov, Maximal B-regular
integro-differential equations, Chin. Ann. Math. Ser. B, 30B(1), (2008),
39-50

\item V. B. Shakhmurov, R.V. Shahmurov, Sectorial operators with convolution
term, Math. Inequal. Appl., V.13 (2), 2010, 387-404.

\item V. B. Shakhmurov, H. K. Musaev, Separability properties of
convolution-differential operator equations in weighted $L_{p}$ spaces,
Appl. and Compt. Math. 14(2) (2015), 221-233.

\item Shakhmurov V. B., Embedding and separable differential operators in
Sobolev-Lions type spaces, Mathematical Notes, (84) 6, (2008), 906-926.

\item Shakhmurov V. B., Maximal regular abstract elliptic equations and
applications, Siberian Math. Journal, (51) 5, (2010), 935-948.

\item Triebel H., Interpolation theory, Function spaces, Differential
operators, North-Holland, Amsterdam, 1978.

\item Shi A., Bai, Yu., Existence and uniqueness of solution to two-point
boundary value for two-sided fractional differential equations, Applied
Mathematics, 4 (6), (2013), 914-918.

\item Weis L, Operator-valued Fourier multiplier theorems and maximal $L_{p}$
regularity, Math. Ann. 319(2001), 735-758.

\item Yakubov S., Yakubov Ya., Differential-operator Equations. Ordinary and
Partial \ Differential Equations, Chapman and Hall /CRC, Boca Raton, 2000.
\end{enumerate}

\end{document}